\documentclass{ifacconf}

\usepackage{natbib}            
\usepackage{graphicx}          
\usepackage{amsmath}
\usepackage{amssymb}
\usepackage{bm}
\usepackage{color}

\begin{document}

\begin{frontmatter}

\title{Path Integral Formulation of Stochastic Optimal Control with Generalized Costs\thanksref{footnoteinfo}} 

\thanks[footnoteinfo]{This work was supported by the NSF CPS project ActionWebs under grant number 0931843, NSF CPS project FORCES under grant number 1239166, by the Director, Office of Science, Computational and Technology Research, U.S. Department of Energy under Contract No. DE-AC02-05CH11231, and by the National Science Foundation under grant DMS-1217065.}

\author[First]{Insoon Yang} 
\author[Second]{Matthias Morzfeld} 
\author[First]{Claire J. Tomlin}
\author[Second]{Alexandre J. Chorin}

\address[First]{Department of Electrical Engineering and Computer Sciences, University of California, Berkeley, CA 94720 USA (e-mail: iyang@eecs.berkeley.edu, tomlin@eecs.berkeley.edu)}                                              
\address[Second]{Department of Mathematics, University of California, Berkeley and Lawrence Berkeley National Laboratory, CA 94720 USA (e-mail: mmo@math.lbl.gov, chorin@math.berkeley.edu)}


\begin{abstract}                          
Path integral control solves a class of stochastic optimal control problems with a Monte Carlo (MC) method for an associated Hamilton-Jacobi-Bellman (HJB) equation. 
The MC approach avoids the need for a global grid of the domain of the HJB equation and, therefore, path integral control is in principle applicable to control problems of moderate to large dimension.
The class of problems path integral control can solve, however, is defined by requirements on the cost function, the noise covariance matrix and the control input matrix. We relax the requirements on the cost function by introducing a new state that represents an augmented running cost. In our new formulation the cost function can contain stochastic integral terms and linear control costs, which are important in applications in engineering, economics and finance. We find an efficient numerical implementation of our grid-free MC approach and demonstrate its performance and usefulness in examples from hierarchical electric load management. The dimension of one of our examples is large enough to make classical grid-based HJB solvers impractical.
\end{abstract}

\end{frontmatter}

\section{Introduction}

Stochastic optimal control is an important tool for controlling or analyzing stochastic systems in many fields, including engineering, economics and finance. 
In electric power systems, for example, uncertainties from renewable energy sources can be managed using stochastic optimal control techniques (e.g.~\cite{Anderson2011, Yang2013}).

A popular method for stochastic optimal control is to compute an optimal feedback map from the system state to the control via dynamic programming (\cite{Bellman1956}). 
For systems modeled by stochastic differential equations (SDEs), dynamic programing leads to a Hamilton-Jacobi-Bellman (HJB) equation, whose viscosity solution corresponds to the value function of the stochastic optimal control problem (\cite{Fleming2005}). An optimal control strategy can be synthesized from information about the value function.

The complexity and memory requirements of grid-based partial differential equation (PDE) solvers increase exponentially as the dimension of the system increases. This 
makes grid-based methods impractical for problems of large dimension. As an alternative to grid-based PDE solvers one can use Monte Carlo (MC) schemes. This is the main idea in ``path integral control'' (\cite{Kappen2005a,Kappen2005b,Todorov2009,Theodorou2010,Theodorou2011,Morzfeld2013}), which use the Feynman-Kac formula (see e.g.~\cite{Chorin2013}), i.e.~a path integral, to solve a class of HJB equations. The current theory of path integral control requires several assumptions about the cost function, noise covariance and control input matrices. These assumptions are often restrictive and can prevent its broad application. 
A related method that can help with the scaling of grid-based PDE solvers is the ``max-plus'' method (\cite{McEneaney2008,McEneaney2009}).

Our goal is to extend path integral control to handle more general cost functions, in particular costs involving stochastic (It\^{o}) integral terms or linear control costs. Cost functions with stochastic integral terms are important in dynamic contract problems where one needs to incorporate the risk aversions of a principal and an agent into the contract (e.g.~\cite{Holmstrom1987}); linear terms in the control cost are important in hierarchical electric load management problems (see Section 4). We introduce a new state variable that integrates the running costs from the initial time to the current time. We then use dynamic programming to derive the HJB equation associated with the control problem in the augmented state space. We can use the Feynman-Kac formula to solve the HJB equation because we do not introduce any additional requirements on the control input matrix or noise covariance matrix. Thus, we can solve control problems with generalized costs using the path integral-based method. 

We use implicit sampling (\cite{Chorin2009,Chorin2010b,Morzfeld2011,Morzfeld2012,Atkins2013,Morzfeld2013}) to evaluate the Feynman-Kac formula numerically. Implicit sampling is an MC sampling scheme that guides the samples to high-probability regions via a numerical optimization. This can be particularly advantageous (compared to traditional MC sampling) if the running cost has deep wells (\cite{Kappen2005a,Morzfeld2013}). 

The efficiency and usefulness of our approach is demonstrated with an example from hierarchical management of groups of thermostatically controlled loads (TCLs). The dimension of the state space of this example is large enough to make grid based PDE solvers impractical.

The remainder of this paper is organized as follows. In Section \ref{review}, we review path integral control and its limitations.
We explain our generalization in Section \ref{path} and also discuss our implementation with implicit sampling. In Section \ref{ex} we apply our method to a problem in direct load control.

\section{Review of Path Integral Control} \label{review}

Consider a process $\{x_t\}_{0\leq t \leq T}$, $x_t \in \mathbb{R}^n$ that is driven by the stochastic differential equation (SDE):
\begin{equation} \label{dynamics1}
\begin{split}
dx_t &= [f(t, x_t) + K u_t ]dt + \sigma dW_t\\
x_0 &= x^0,
\end{split}
\end{equation}
where $K \in \mathbb{R}^{n \times m}$ and $\sigma \in \mathbb{R}^{n \times n}$, and $\{u_t\}_{0 \leq t \leq T}$, $u_t \in \mathbb{R}^m$ is the control process; $\{W_t\}_{t \geq 0}$ is an $n$-dimensional Brownian motion on a suitable probability space. 

Path integral control requires a cost function of the form 
\begin{equation} \label{originalcost}
\min_{u \in \mathbb{U}} \quad \mathbb{E} \left [ \Phi(x_T) +  \int_0^T V(t, x_t)dt +  \frac{1}{2}\int_0^T u_t^\top R u_t \: dt \right ]
\end{equation}
where $\mathbb{U}$ is the set of feasible controls, $R\in \mathbb{R}^{m\times m}$ is a given positive definite matrix and $\Phi(x_T)$ and $V(t,x_t)$ are the terminal and running costs, respectively (\cite{Kappen2005a, Theodorou2010}); the expected value is taken over trajectories of \eqref{dynamics1}.

Path integral control further requires that there exists a constant $\lambda$ such that
\begin{equation} \label{assumption}
\sigma \sigma^\top = \lambda KR^{-1}K^\top.
\end{equation}
This condition implies that the control input through a subsystem (or a channel) with a higher noise variance is cheaper than that through a subsystem with a lower noise variance (\cite{Kappen2005a}).

Dynamic programming can be used to derive the HJB equation for a stochastic optimal control problem. This HJB equation is generally nonlinear, however for the class of problems defined by \eqref{originalcost} and \eqref{assumption}, a log-transformation can be used to derive a linear HJB equation (\cite{Kappen2005a}). The solution of this linear HJB equation can be found via the Feynman-Kac formula (a path integral). 

In computations, one discretizes the path to obtain a convergent approximation of the Feynman-Kac formula. MC sampling can be used to evaluate this discretized path integral numerically and, thus, to find the solution of the HJB equation. Specifically, MC schemes  evaluate the solution of the HJB equation locally without requiring the solution nearby, so that there is no need for a (global) grid of the domain. This feature makes path integral control with MC sampling particularly effective and memory efficient when applicable.

\section{Path Integral Formulation of Stochastic Optimal Control with Generalized Costs} \label{path}

Our goal is to develop a path integral method that can handle cost functions that are more general than \eqref{originalcost}. 
Specifically, consider a stochastic optimal control problem of the form
\begin{equation}\label{stoch_opt1}
\begin{split}
\min_{u \in \mathbb{U}} \quad &\mathbb{E} \left [ \Phi \left (x_T, \int_0^T V(t, x_t)dt + \int_0^T c(t) dx_t \right ) \right.\\
&\left. \quad + \frac{1}{2}\int_0^T u_t^\top R u_t \: dt \right ],
\end{split}
\end{equation}
where $\Phi : \mathbb{R}^{n} \times \mathbb{R} \to \mathbb{R}$ is a given function, and where $V(t,x_t) \in \mathbb{R}$ and $c(t) \in \mathbb{R}^{1 \times n}$ are the running costs, the latter being the running cost associated with a stochastic state trajectory. This generalized cost function has been used in dynamic incentive design problems to take into account the principal's and the agent's risk aversions (\cite{Holmstrom1987}, \cite{Ou-Yang2003}). We will further show that the generalized cost function can include running costs that are linear in the control and discuss its usefulness in hierarchical management for thermostatically controlled loads (in Section \ref{ex}).

\subsection{The augmented system}

We reformulate the stochastic optimal control problem with generalized cost \eqref{stoch_opt1} to a related problem for which we can derive the associated HJB equation. The key idea is to introduce the new state variable
\begin{equation}  \label{new}
y_t := \int_0^t V(s, x_s) ds + \int_0^t c(s) dx_s,
\end{equation}
where the stochastic integral is in the It\^{o} sense. Then, the process $\{ y_t \}_{0 \leq t \leq T}$ is driven by the following SDE:
\begin{equation} \label{dynamics2}
\begin{split}
d y_t &= [V(t, x_t) + c(t)( f(t, x_t) +K u_t) ] dt + c(t) \sigma dW_t\\
y_0 &= 0,
\end{split}
\end{equation}
where $\{W_t\}_{t \geq 0}$ is the Brownian motion in \eqref{dynamics1}. 
We can define an augmented state space by $(x_t,y_t)$. 
The stochastic optimal control problem \eqref{stoch_opt1} can be rewritten in this augmented state space as
\begin{equation}\label{stoch_opt2}
\min_{u \in \mathbb{U}} \quad \mathbb{E} \left [ \Phi \left (x_T, y_T \right )  + \frac{1}{2}\int_0^T u_t^\top R u_t \:dt \right ]
\end{equation}
where the expectation is over the trajectories of \eqref{dynamics1} and \eqref{dynamics2}.

\subsection{The HJB equation of the augmented system}
We will now provide the details about the derivation of the HJB equation for the augmented state space problem. This will highlight that path integral control can be applied to problems with generalized costs \eqref{stoch_opt1} provided that \eqref{assumption} holds, while no additional restrictions are required.

We define the value function of the stochastic optimal control as
\begin{equation} \label{value}
\begin{split}
&\phi(\bm{x}, \bm{y}, t) \\
& := \min_{u \in \mathbb{U}} \: \mathbb{E}_{\bm{x}, \bm{y}, t} \left [  \Phi(x_T, y_T) + \frac{1}{2} \int_t^T u_s^\top R u_s ds \right ], 
\end{split}
\end{equation}
where the expected values is taken over trajectories of \eqref{dynamics1} and \eqref{dynamics2} starting at $(x_t, y_t) = (\bm{x}, \bm{y})$.
Let
\begin{equation}  \nonumber
 \Lambda (t) := \begin{bmatrix}
\bold{I}_{n} \\
c(t) 
\end{bmatrix} \in \mathbb{R}^{(n+1) \times n},
\end{equation}
where $\bold{I}_n$ denotes the $n \times n$ identity matrix,
and $\hat{\sigma}(t) := \Lambda (t) \sigma$.
Dynamic programming gives the HJB equation
\begin{equation} \label{HJB1}
\begin{split}
& \phi_t + \min_{a \in \mathbb{R}^{m}} \left \{ \frac{1}{2} \mbox{tr}(\hat{\sigma} \hat{\sigma}^\top D^2 \phi) + (f(t, \bm{x}) + K a)^\top D_{\bm{x}} \phi \right.\\ 
& \;\;\; \left. + (V(t, \bm{x}) + c(t) f(t,\bm{x}) +  c(t)K a) D_{\bm{y}} \phi  +\frac{1}{2} a^\top R a  \right \} = 0\\
&\phi(\bm{x}, \bm{y},T) = \Phi(\bm{x}, \bm{y}).
\end{split}
\end{equation}
The minimizer in \eqref{HJB1} is given by
\begin{equation} \label{optControl}
a^* = -R^{-1} \hat{K}(t)^\top D \phi,
\end{equation}
where $\hat{K} (t) = \Lambda (t) K$ and $D\phi := (D_{\bm{x}} \phi, D_{\bm{y}} \phi) \in \mathbb{R}^{n+1}$. Plugging this minimizer into \eqref{HJB1}, we obtain
\begin{equation} \label{HJB2}
\begin{split}
&\phi_t + \frac{1}{2} \mbox{tr} (\hat{\sigma} \hat{\sigma}^\top D^2 \phi) + f(\bm{x})^\top D_{\bm{x}}\phi + (V(\bm{x}) + c f(\bm{x}) )^\top D_{\bm{y}} \phi\\ 
&- \frac{1}{2} D \phi^\top \hat{K} R^{-1} \hat{K}^\top D \phi = 0,
\end{split}
\end{equation}
with the same terminal condition.

Following \cite{Kappen2005a}, we define a new value function by
\begin{equation}\nonumber
\psi(\bm{x},\bm{y}, t) = \exp \left (-\frac{\phi(\bm{x},\bm{y},t)}{\lambda} \right ),
\end{equation}
for all $(\bm{x},\bm{y}, t) \in \mathbb{R}^{n}\times \mathbb{R} \times [0,T]$,
where the constant $\lambda$ satisfies the condition \eqref{assumption}.
A calculation shows that \eqref{assumption} implies 
\begin{equation} \label{condition2}
\begin{split}
\hat{\sigma} (t) \hat{\sigma} (t)^\top  &= \Lambda(t) \sigma \sigma^\top \Lambda(t) = \lambda \hat{K}(t) R^{-1} \hat{K}(t)^\top.
\end{split}
\end{equation}
Substituting $\phi = -\lambda \log \psi$ in \eqref{HJB2} and using \eqref{condition2} gives a linear PDE for $\psi$: 
\begin{equation} \label{HJB3}
\begin{split}
&\psi_t + \frac{1}{2} \mbox{tr} (\hat{\sigma} \hat{\sigma}^\top D^2 \psi) 
 + f(\bm{x})^\top D_{\bm{x}} \psi \\
& \qquad \qquad \qquad \quad \quad \: + (V(\bm{x}) + c f(\bm{x}))^\top D_{\bm{y}} \psi
 =0\\ 
&\psi (\bm{x}, \bm{y},T) = \exp \left ( -\frac{1}{\lambda} \Phi(\bm{x}, \bm{y}) \right ).
\end{split}
\end{equation} What we have shown is that  the condition \eqref{assumption} is sufficient to obtain a linear HJB equation \eqref{HJB2} for the case of generalized costs, i.e.~path integral control is applicable to problems with generalized cost functions provided that \eqref{assumption} holds.

\subsection{Implementation with implicit sampling} \label{genpath}
Upon solving \eqref{HJB2} or \eqref{HJB3} locally around $(x_t, y_t) = (\bm{x}, \bm{y})$ at time $t$, an optimal control at time $t$ can be determined by 
\begin{equation} \nonumber
u_t^* = -R^{-1} \hat{K}^\top D\phi (\bm{x}, \bm{y}, t).
\end{equation}

Instead of solving \eqref{HJB3} with a grid-based scheme, we use the Feynman-Kac formula (e.g.~\cite{Chorin2013}):
\begin{equation} \label{exp}
\psi(\bm{x},\bm{y}, t) = \mathbb{E} \left [ \exp  \left( -\frac{1}{\lambda} \Phi(\tilde{x}_T, \tilde{y}_T)\right )  \right ],
\end{equation}
where the expectation is taken over the stochastic processes $\{\tilde{x}_s\}_{0 \leq s \leq T}$ and $\{\tilde{y}_s\}_{0 \leq s \leq T}$ driven by
\begin{equation} \label{stochd}
\begin{split}
d\tilde{x}_s &= f(s, \tilde{x}_s) ds + \sigma dW_s \\
d\tilde{y}_s &= [V(s, \tilde{x}_s) + c(s) f(s, \tilde{x}_s)] ds + c(s) \sigma dW_s\\
(\tilde{x}_t, \tilde{y}_t)&= (\bm{x}, \bm{y}).
\end{split}
\end{equation}
Note that the Brownian motion in the dynamics of $\tilde{y}$ is the same as that in the dynamics of $\tilde{x}$.
The Feynman-Kac formula holds because the HJB equation \eqref{HJB3} is linear and parabolic. We use \eqref{exp} to locally evaluate $\psi$ at any $(\bm{x}, \bm{y},t)$ without a grid (which is impractical in high dimensions). 

We first approximate the Feynman-Kac formula, and discretize the time interval $[0,T]$ with $\{t_i\}_{i=0}^M$, where
$t_i = i \Delta t$ and $\Delta t : = \frac{T}{M}$.
Let $\{\bold{x}_i\}_{i=j}^M$ and $\{\bold{y}_i\}_{i=j}^M$ be the discretized trajectories \eqref{stochd} of $\tilde{x}_s$ and $\tilde{y}_s$ for $s \in [t_j,T]$, respectively, starting from $(\bold{x}_j, \bold{y}_j) = (\bm{x}, \bm{y})$.
Recall that the Brownian motion in the dynamics of $\tilde{y}$ is the same as that in the dynamics of $\tilde{x}$. Therefore, the uncertainty in the path $\{\bold{y}_i\}_{i=j}^M$ is determined by the uncertainty in the path $\{\bold{x}_i\}_{i=j}^M$ and the probability distribution of the path $\{\bold{x}_i\}_{i=j}^M$ defines the probability of the paths $\{\bold{x}_i\}_{i=j}^M$ and $\{\bold{y}_i\}_{i=j}^M$.
We can thus use
\begin{equation}\label{pi}
\begin{split}
\psi(\bm{x},\bm{y},t_j) &\approx \int d\bold{x}_{j+1} \dots \int d\bold{x}_M \: p(\bold{x}_{j+1}, \dots, \bold{x}_M) \\
& \times \exp \left ( -\frac{1}{\lambda} \Phi (\bold{x}_M, \bold{y}_M)  \right ).
\end{split}
\end{equation}
to approximate \eqref{exp}. 
When computing $\bold{y}_M$, we can use the fact that $\tilde{y}$ satisfies \eqref{stochd} to find that
\begin{equation} \label{yy}
\bold{y}_M = \bold{y}_j + \sum_{i = j+1}^M V(t_i, \bold{x}_i) \Delta t + \sum_{i = j+1}^M c(t_i) (\bold{x}_i - \bold{x}_{i-1}),
\end{equation}
converges to $\tilde{y}_T$ of the second equation in \eqref{stochd} as $\Delta t \to 0$.

To evaluate the discretized path integral \eqref{pi} we use MC sampling. For example, we can generate $Q$ trajectories $\{  (\bold{x}_i^{(q)}, \bold{y}_i^{(q)}  ) \}_{i=j}^M$ of \eqref{stochd} starting from $(\bold{x}_j^{(q)}, \bold{y}_j^{(q)}) = (\bm{x}, \bm{y})$, $q= 1, \dots, Q$. 
 The MC approximation of \eqref{pi} then becomes 
\begin{equation} \label{MCapprox}
\psi(\bm{x}, \bm{y},t_j) \approx \frac{1}{Q}\sum_{q=1}^Q \exp \left ( -\frac{1}{\lambda} \Phi \left (\bold{x}_M^{(q)}, \bold{y}_M^{(q)} \right )  \right ).
\end{equation}
As the number of samples $Q$ goes to infinity, the right-hand side converges weakly to the value function, and an optimal control can be synthesized as
\begin{equation} \nonumber
u_{t_j}^* = -R^{-1} \hat{K}(t_j)^\top  D \phi (\bm{x},\bm{y}, t_j),
\end{equation}
where $D \phi (\bm{x},\bm{y}, t_j)$ is computed via numerical differentiation (e.g.~forward or centered differences). In practice, we cannot use ``infinitely'' many samples. In fact we may only be able to generate a few ($O(1000)$) samples. In this case, this ``standard MC scheme" can fail, in particular if the running cost $b$ has deep wells (\cite{Kappen2005a, Theodorou2010,Morzfeld2013}). Deep but thin wells in $V$ imply that many trajectories of \eqref{stochd} end up where $V$ is large and thus contribute little to the approximation \eqref{MCapprox} of $\psi$. This pitfall can be overcome if one guides the samples to remain where the running cost $V$ is small.

Here, we achieve this guiding effect via implicit sampling, which takes the running cost into account when generating the samples (\cite{Chorin2009,Chorin2010b,Morzfeld2011,Morzfeld2012,Atkins2013,Morzfeld2013}). 
To use implicit sampling, we first substitute $p(\bold{x}_{j+1}, \dots, \bold{x}_M)$ in the path integral formula \eqref{pi} 
\begin{equation} \nonumber
\begin{split}
p(\bold{x}_{j+1}, \dots, \bold{x}_M) &= \prod_{i = j+1}^M p(\Delta \bold{x}_i)\\
&= (2\pi \mbox{det}(\sigma \sigma^\top \Delta t))^{-(M-j)/2} \\
& \times \exp \left (-\frac{1}{2} \sum_{i=j+1}^M \bold{q}_i^\top (\sigma \sigma^\top \Delta t)^{-1} \bold{q}_i \right ),
\end{split}
\end{equation}
where $\bold{q}_i := \Delta \bold{x}_i - f(t_{i-1}, \bold{x}_{i-1}) \Delta t $ and
\begin{equation} \label{nonumber}
\Delta \bold{x}_i := \bold{x}_i - \bold{x}_{i-1} \sim \mathcal{N}( f(t_{i-1}, \bold{x}_{i-1}) \Delta t, \sigma \sigma^\top \Delta t).
\end{equation}
We can then rewrite \eqref{pi} as
\begin{equation} \nonumber
\begin{split}
\psi(\bm{x}, \bm{y}, t_j) &\approx (2\pi \mbox{det}(\sigma \sigma^\top \Delta t))^{-(M-j)/2} \\
&\quad \int d\bold{x}_{j+1} \cdots \int d\bold{x}_M \exp (-F (\bold{x})),
\end{split}
\end{equation}
where $\bold{x}:= (\bold{x}_{j+1}, \dots, \bold{x}_M) \in \mathbb{R}^{n (M-j)}$ and
\begin{equation} \nonumber
\begin{split}
&F(\bold{x}) = \frac{1}{\lambda} \Phi (\bold{x}_M, \bold{y}_M) +\frac{1}{2} \sum_{i=j+1}^M \left [ \bold{q}_i^\top  (\sigma \sigma^\top \Delta t)^{-1} \bold{q}_i \right ]
\end{split}
\end{equation}
and $\bold{y}_M$ is given by \eqref{yy}. Note that the running cost $V$ is included in $F$, i.e.~$F$ contains information about where $V$ is large and where it is small.

The high-probability samples (those that correspond to a small running cost) are in the neighborhood of where $F$ is small, i.e.~in the neighborhood of the minimizer
\begin{equation} \nonumber
\mu:= \arg\min_{\bold{x}} \: F(\bold{x}).
\end{equation}
We can generate samples in this neighborhood by mapping the high probability region of a ``reference variable'' $\xi$ to this neighborhood. Here, we choose a Gaussian reference variable $\xi\sim\mathcal{N}(0, \bold{I}_{n(M-j)})$, and a linear mapping 
\begin{equation} \label{mechanism}
\bold{x} = \mu + L^{-1} \xi,
\end{equation}
where $L$ is a Cholesky factor of the Hessian $H=LL^\top$ of $F$ evaluated at the minimizer $\mu$. 

To see why \eqref{mechanism} maps the high probability region of $\xi$ to the neighborhood of $\mu$, note that \eqref{mechanism} implies that 
\begin{equation}\label{approxSamplingEq}
\tilde{F}(\bold{x}) - \zeta = \frac{1}{2} \xi^\top \xi,
\end{equation}
where $\zeta := \min_{\bold{x}} \: F(\bold{x})$ and where
\begin{equation} \nonumber
\tilde{F}(\bold{x}) = \zeta + \frac{1}{2} (\bold{x} -\mu)^\top H (\bold{x} - \mu)
\end{equation}
is the Taylor expansion to order $2$ of $F$. High-probability samples of $\xi$ are in the neighborhood of the origin, and hence, the right hand side of \eqref{approxSamplingEq} is small with  high probability. Accordingly, the left hand side of \eqref{approxSamplingEq} is also small with  high probability and, therefore, the solutions of \eqref{approxSamplingEq} are close to the minimizer $\mu$ of $\tilde{F}$, which is also the minimizer of $F$.

The sampling equation \eqref{mechanism} defines an invertible change of variables from $\bold{x}$ to $\xi$, which can be used to write \eqref{pi} in terms of the reference variable $\xi$ 
\begin{equation} \nonumber
\begin{split}
&\psi(\bm{x}, \bm{y}, t_j) \approx (2\pi \mbox{det}(\sigma \sigma^\top \Delta t))^{-\frac{M-j}{2}}  \int d\bold{\xi}_{1} \cdots \int d\bold{\xi}_{n(M-j)}\\
&\qquad \quad \times \exp \left (-\zeta + \tilde{F}(\bold{x}(\xi)) - F(\bold{x}(\xi)) - \frac{1}{2} \xi^\top \xi \right ) \mbox{det}(J),
\end{split}
\end{equation}
where $J = {\partial \bold{x}}/{\partial \xi}=L^{-1}$. Since $\xi \sim \mathcal{N}(0, \bold{I}_{n(M-j)})$, we may also write
\begin{equation} \nonumber
\psi(\bm{x}, \bm{y}, t_j) \approx  \frac{\mathbb{E}_{\xi} [ \exp (-\zeta +  \tilde{F}(\bold{x}) - F(\bold{x}))]}{(\mbox{det}(\sigma \sigma^\top \Delta t))^{\frac{M-j}{2}}\mbox{det}(L)}.
\end{equation}
Our implicit sampling-based path integral algorithm can be summarized in three steps:
\begin{enumerate}
\item
Generate samples, $\xi^{(q)}$, $q=1, \dots, Q$, from $\xi \sim \mathcal{N}(0, \bold{I}_{n(M-j)})$;
\item
Transform the samples through the mechanism \eqref{mechanism}, i.e., $\bold{x}^{(q)} = \mu + L^{-1} \xi^{(q)}$;
\item
Evaluate the value function $\psi$ at $(\bm{x}, \bm{y}, t_j)$ as
\begin{equation}\nonumber
\psi (\bm{x}, \bm{y}, t_j) \approx \frac{1}{Q} \sum_{q = 1}^Q  \frac{\exp(-\zeta + \tilde{F}(\bold{x}^{(q)}) - {F}(\bold{x}^{(q)}))}{(\mbox{det}(\sigma \sigma^\top \Delta t))^{\frac{M-j}{2}}\mbox{det}(L)}.
\end{equation}
\end{enumerate}
For more detail about implicit sampling in the context of path integral control see \cite{Morzfeld2013}.

\section{Hierarchical Management for Stochastic Thermostatically Controlled Loads} \label{ex}

We consider a scenario in which a utility company manages the operation of thermostatically controlled loads (TCLs), such as air conditioners, to minimize payments in the real-time market and the discomfort of the consumers.
If the utility controls all of the TCLs by monitoring all of the system states (e.g.~indoor temperatures), the synthesis of an optimal control is intractable due to high system dimension. 

\begin{figure}[tb] 
\begin{center}
\includegraphics[width = 2.9in]{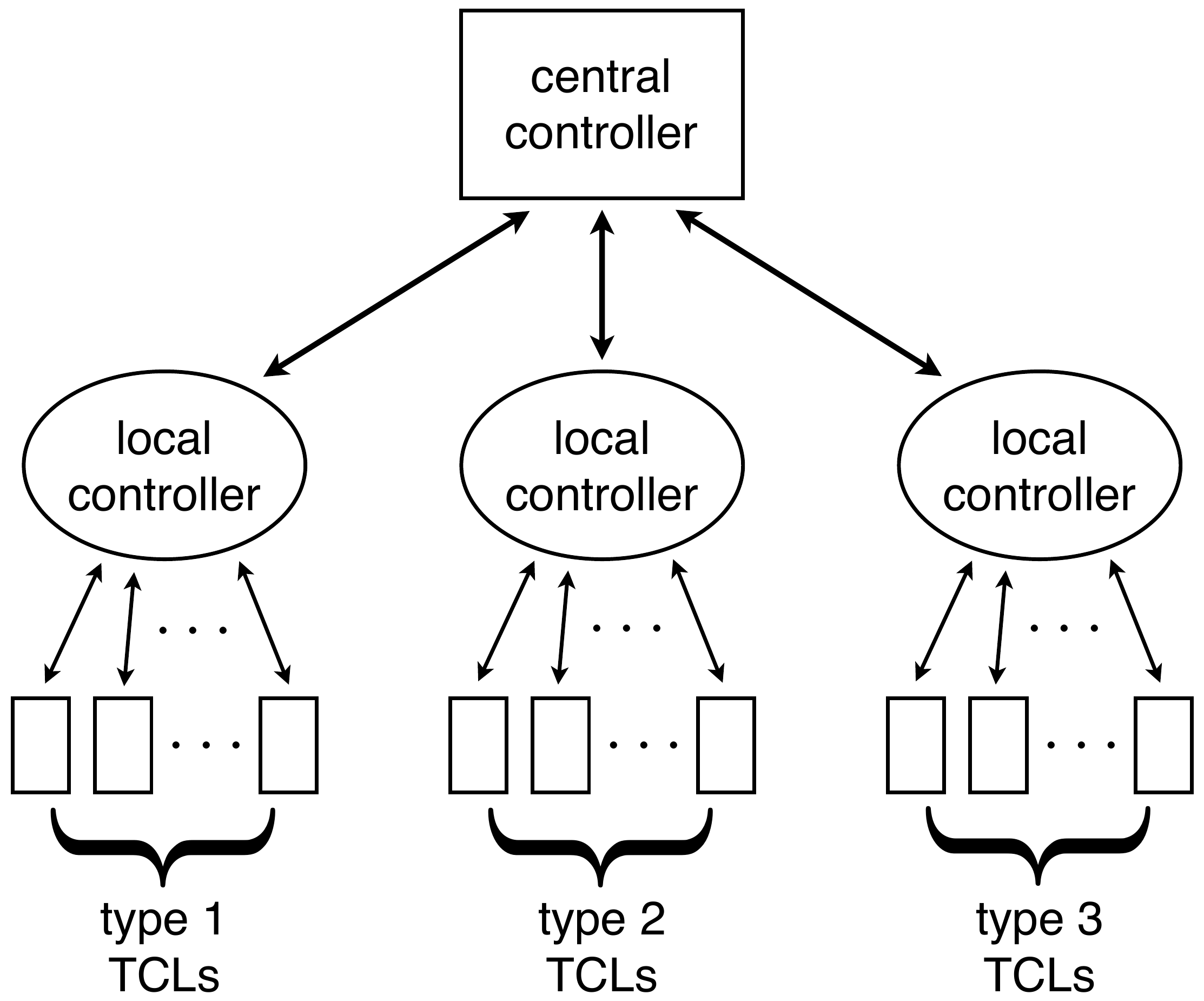}
\caption{Hierarchical control for three types of TCLs.
}
 \label{fig:hierarchical}
 \end{center}
\end{figure}

To overcome the system dimension issue,
we propose the hierarchical management framework in Fig. \ref{fig:hierarchical}. A similar hierarchical load control framework is proposed in \cite{Callaway2011}.
We assume that all of the TCLs are air conditioners and that exactly one TCL controls the indoor temperature of a single room.
The TCLs are classified into $n$ types, which differ in their temperature-TCL model characteristics (see section \ref{tclmodel}).
The key components of the hierarchical management are the local controllers. Local controller $i$ determines which subset of type $i$ TCLs should be turned on at each time $t \in [0, T]$. Local controller $i$ communicates with type $i$ TCLs to gather the data of the indoor temperatures, and then ranks the TCLs from the highest to lowest discomfort level (the precise definition of the discomfort level is given Section \ref{controlobjectives}). Each local controller also reports the average indoor temperature to the central controller.

The central controller, which communicates only with the local controllers, 
computes each TCL group's optimal power consumption by considering the average temperature-TCL models of all of the groups.
Note that the total dimension of the averaged systems is much less than that of all of the temperature-TCL systems. The central controller then delivers group $i$'s optimal power consumption to local controller $i$ for $i = 1, \dots, n$.

After receiving the optimal power consumption level, local controller $i$ determines the number $\mathcal{O}^*$ of TCLs that should be turned on by dividing the optimal level by one TCL's nominal power consumption. 
The local controller changes the operation of the TCL with the highest rank (the one with the highest discomfort level): i.e., if the TCL is off, the local controller turns it on. The local controller then considers the TCL with the second highest rank, and so on. It stops making changes once the number of TCLs that are on is equal to $\mathcal{O}^*$.

We consider the central controller's problem of determining the optimal power consumptions of 
all of the TCL groups. This stochastic optimal control problem is solved by the generalized path integral method with implicit sampling.

\subsection{Stochastic indoor temperature-TCL Model} \label{tclmodel}

\begin{figure}[tb] 
\begin{center}
\includegraphics[width = 3.5in]{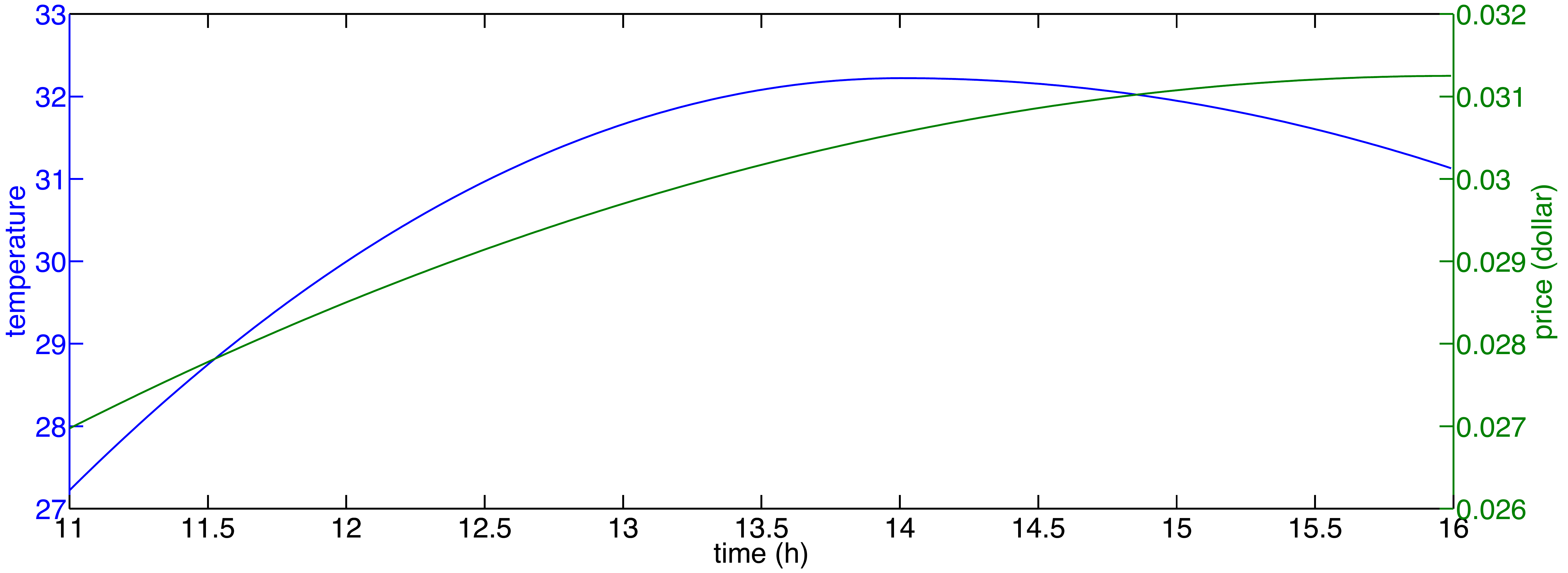}
\caption{A forecast of outdoor air temperature, $\Theta_t^O$ (blue), and an energy price profile, $p(t)$, in the real-time market (red).
}
 \label{fig:outdoor}
 \end{center}
\end{figure}

Let $\Theta_t^{i,j}$ be the indoor air temperatures at time $t$ that are controlled by TCL $j$ of type $i$  for $j = 1, \dots, N_i$, where $N_i$ is the number of type $i$ TCLs.
The equivalent thermal parameter (ETP) model suggests that the dynamics of the temperatures depend on the conductance $R_{i,j}$ between the outdoor air and indoor air and the conductance $C_{i,j}$ of the indoor air and the thermal mass (\cite{Sonderegger1978}).
We assume that TCLs of the same type have the same thermal parameters: i.e., $R_i := R_{i,j}$ and $C_i := C_{i,j}$  for $j = 1, \dots, N_i$.
Let $\Theta_t^{i}$ be the average indoor air  temperature at time $t$, i.e., $\Theta_t^{i} := \frac{1}{N_i} \sum_{j=1}^{N_i} \Theta_t^{i,j}$.
We consider the following stochastic model for the average temperatures: 
\begin{equation}\label{tcl}
\begin{split}
d\Theta_t^{i} &= \left [ \alpha_i(\Theta_t^{O,i} - \Theta_t^{i})  -  \kappa_i \frac{u_t^i}{N_i} \right] dt + \sigma_i dW_t^i\\
\end{split}
\end{equation}
where $\alpha_i = R_i/ C_i$.
Here, $u_t^i$ denotes the total power consumption of the type $i$ TCLs
and
$\Theta_t^{O,i}$ denotes the (day-ahead) forecast of the average outdoor air temperature at time $t$. We use the forecast shown in Fig. \ref{fig:outdoor} for our simulations.
The effect of the forecast error is modeled by the stochastic term $\sigma_i dW_t^i$.

Let $x_t^{i} :=  \Theta_t^{i}$ for $i = 1, \dots, n$, $u_t := (u_t^1, \dots, u_t^n)$, $W_t := (W_t^1, \dots, W_t^n)$, and let $K$ and $\sigma$ be $n \times n$ diagonal matrices whose $i$th diagonal entries are $-\kappa_i/N_i$ and $\sigma_i$, respectively.
With this notation, we can compactly write the stochastic temperature-TCL model for all types as
\begin{equation} \label{tcl1}
dx_t = [f(t, x_t) +  K u_t] dt + \sigma dW_t,
\end{equation}
where
\begin{equation} \nonumber
\begin{split}
f_{i}(t, \bm{x}) &:=  \alpha_i (\Theta_t^{O,i} - \bm{x}^i)
\end{split}
\end{equation}
for $i = 1, \dots, n$.
The central controller determines the total power consumption vector $u$ using this stochastic temperature-TCL model to achieve the objective discussed below.

\subsection{Control objectives} \label{controlobjectives}

A control objective is to keep the discomfort level small. The discomfort level should be small if the average indoor temperatures are within a given temperature range, $[\underline{\Theta}, \overline{\Theta}]$. The discomfort level increases significantly as the indoor temperature drops below $\underline{\Theta}$ or increases above $\overline{\Theta}$. 
To model this characteristic of the discomfort level, we define
\begin{equation} \nonumber
\bar{b}(x_t) := \eta_1\left (\sum_{i=1}^n  b_i (x_t) \right )^2,
\end{equation}
where
\begin{equation} \nonumber
b_i(x_t) :=  \left [ \exp(\eta_2 (x_t^i - \overline{\Theta})) + \exp(\eta_2 (\underline{\Theta} - x_t^i))  \right ] 
\end{equation}
for some positive constants $\eta_1$ and $\eta_2$.

Another control objective is to minimize the cost for real-time balancing of the supply and demand. Suppose that the day-ahead energy purchases are sufficiently small and therefore any power supplied to the TCLs must be purchased in the real-time market.
The balancing cost in the real-time market at time $t$ is computed as
$p(t) \bold{1} u_t$, 
where $p(t)$ is the unit price of power in kW at time $t$ (see Fig.~\ref{fig:outdoor}) in the real-time market and $\bold{1}$ is a $1$ by $n$ vector whose elements are all $1$'s.

The central controller can determine an optimal power consumption strategy for each TCL group by solving the stochastic optimal control problem:
\begin{equation} \label{opt1}
\begin{split}
\min_{u \in \mathbb{U}} \quad &\mathbb{E} \left [ \int_0^T \bar{b} (x_t) + p(t) \bold{1} u_t
 + \frac{1}{2} u_t^\top R u_t  \:dt  \right ]\\
 \mbox{subject to} \quad &\mbox{\eqref{tcl1}}.
 \end{split}
\end{equation}
Note that this problem cannot be solved with the path integral method as described in Section \ref{review} due to the term $p(t) \bold{1} u_t$, which is linear in control.
However, we can solve this problem with our generalized formulation of path integral control.

Specifically, we use the stochastic temperature-TCL model \eqref{tcl1} and the fact that $\mathbb{E} [ \int_0^T p(t) \bold{1} K^{-1} \sigma dW_t ] = 0$ to rewrite the objective of \eqref{opt1} as
\begin{equation} \nonumber
\begin{split}
\min_{u \in \mathbb{U}} \quad &\mathbb{E} \left [ \int_0^T \bar{b} (x_t) - p(t) \bold{1} K^{-1} f(t, x_t) dt \right. \\
&\left. \quad + \int_0^T  p(t) \bold{1} K^{-1} dx_t
 +\frac{1}{2} \int_0^T  u_t^\top R u_t  \:dt  \right ].
 \end{split}
\end{equation}
Let $V(t, \bm{x}) := \bar{b} (\bm{x}) - p(t) \bold{1} K^{-1} f(t, \bm{x})$ and $c(t) := p(t) \bold{1} K^{-1}$. 
We introduce a new state variable driven by the following SDE:
\begin{equation} \label{tcl2}
\begin{split}
dy_t &= [V(t, x_t) + c(t) (f(t, x_t) + K u_t)] dt + c(t) \sigma dW_t \\
&= [\bar{b}(x_t)  + c(t) K u_t] dt + c(t) \sigma dW_t 
\end{split}
\end{equation}
with $y_0 = 0$.
In augmented state space, we obtain a stochastic optimal control problem with cost function
\begin{equation} \nonumber
\begin{split}
\min_{u \in \mathbb{U}} \quad &\mathbb{E} \left [ y_T
 + \frac{1}{2} \int_0^T u_t^\top R u_t  \:dt  \right ]\\
\mbox{subject to} \quad &\mbox{\eqref{tcl1} and \eqref{tcl2}}.
\end{split}
\end{equation}
We apply the generalized path integral method described in Section \ref{genpath} to solve this problem.

\subsection{Numerical tests}

\begin{figure}[tb] 
\begin{center}
\includegraphics[width = 3.25in]{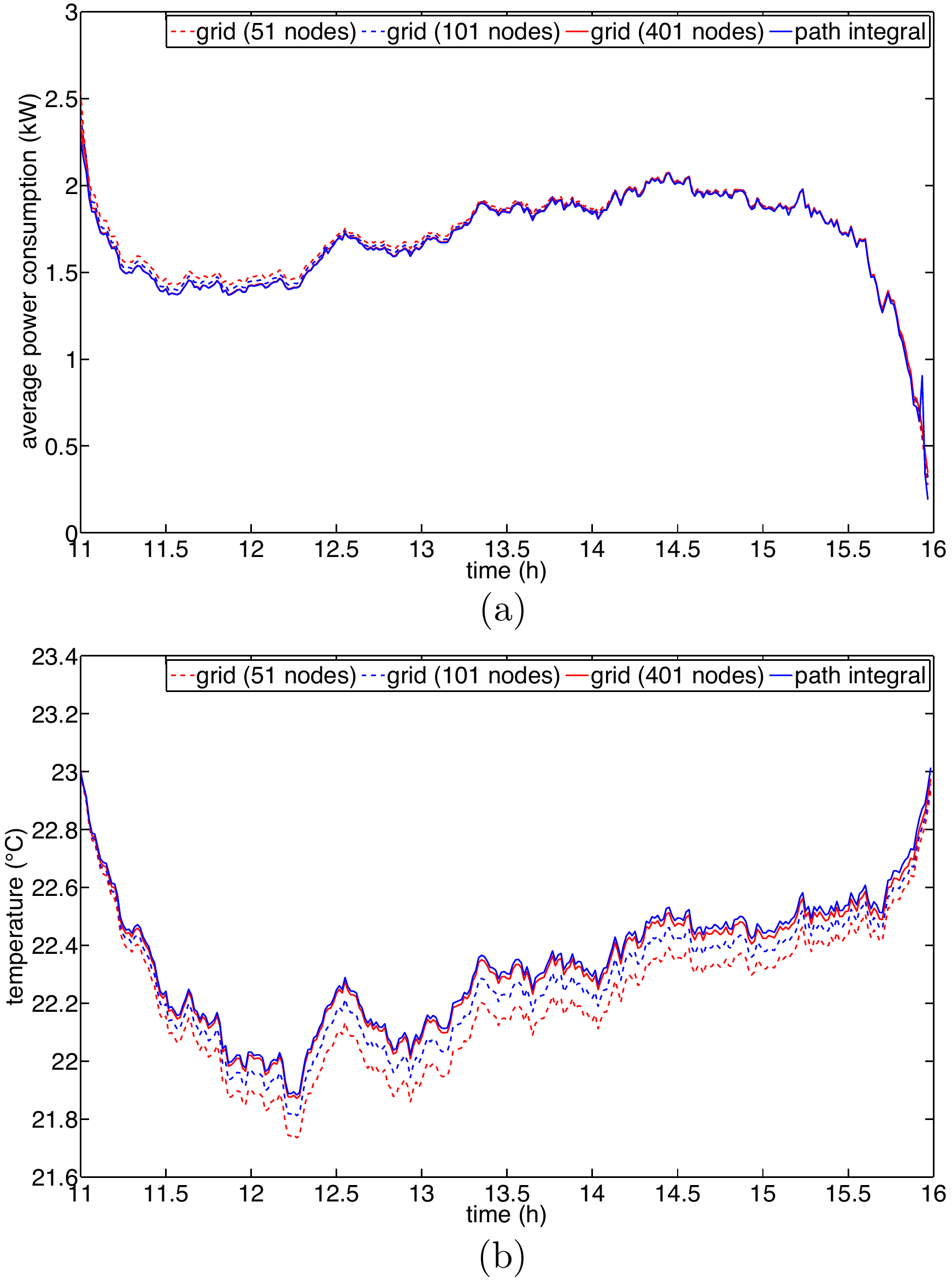}
\caption{Comparison of the grid-based method and the path integral-based method:
(a) computed optimal control and (b) corresponding state trajectory.
}
 \label{fig:one}
 \end{center}
\end{figure}  

\subsubsection{Single TCL type}

We first consider the case in which all of the TCLs are of the same type. The grid-based method is applicable to this simple example so that we can compare results obtained via the grid-based method to path integral control. We set $N_1 = 100$, $\alpha_1 = 0.4834$, $\kappa_1 = 2.5$, $\sigma_1 = 0.1$, $\eta_1 = 2$, $\eta_2 = 5$, $R = 4.6875 \times 10^{-2}$, $(\underline{\Theta}, \overline{\Theta}) = (19,22)$ and $\Theta_0 = 23$.

For the path integral method, we use $300$ nodes for the time discretization of the time interval $[11\mbox{h}, 16\mbox{h}]$, and $5$ samples generated via implicit sampling.
For the grid-based method, we use a varying number of nodes ($51$, $101$ and $401$) to discretize the state space. The varying number of nodes is used to study the convergence of the grid-based method as the discretization is refined.
Both methods lead to a control that induces overcooling before $12$pm when the energy price is low even though the outdoor temperature is not high during this period.

Fig. \ref{fig:one} indicates that the optimal control and corresponding state trajectory obtained by the grid-based method approach those obtained by the path integral method as the number of nodes increases. In fact, The difference between the indoor temperature with the optimal control obtained by the path integral method and that with the optimal control obtained by the grid method (with 401 nodes) is less than $1\%$.
This observation suggests that the path integral method is accurate, i.e., the local solution of the HJB equation we obtain via implicit sampling is the desired viscosity solution because the grid-based solution converges to the viscosity solution as the discretization is refined (e.g.~\cite{Sethian1999,Osher2002} and the references therein).

\subsubsection{Six types of TCLs}

We now consider a more realistic scenario with six types of TCLs. The dimension of this problem is large enough to make the grid-based approach impractical. We set $\alpha = (0.4834, 0.6043, 0.7251,$ $0.8460, 0.9669, 1.0877)$, $\eta_1 = 0.1/6^2$ and $(\underline{\Theta}, \overline{\Theta}) = (18,21.5)$; all other parameters are as above. We use $100$ nodes for the time discretization and $5^6 = 15625$ samples generated by implicit sampling because using $5$ samples gives a good result in the case of single TCL type.

We observe that as $\alpha_i$ increases, the heat conduction from the outside to the inside increases. Because we assume that the outdoor temperature is higher than the indoor temperature (see Fig. \ref{fig:outdoor}), a TCL in a room with a higher $\alpha_i$ consumes more power to reduce the discomfort level. As expected, the optimal power consumption of the type $6$ TCLs (those with the highest $\alpha$) is the highest, and that of the type $1$ TCLs (those with the lowest $\alpha$) is the lowest (see Fig. \ref{fig:one} (a)). This is in line with what we intuitively expect from an optimal control.

\begin{figure}[tb] 
\begin{center}
\includegraphics[width = 3.25in]{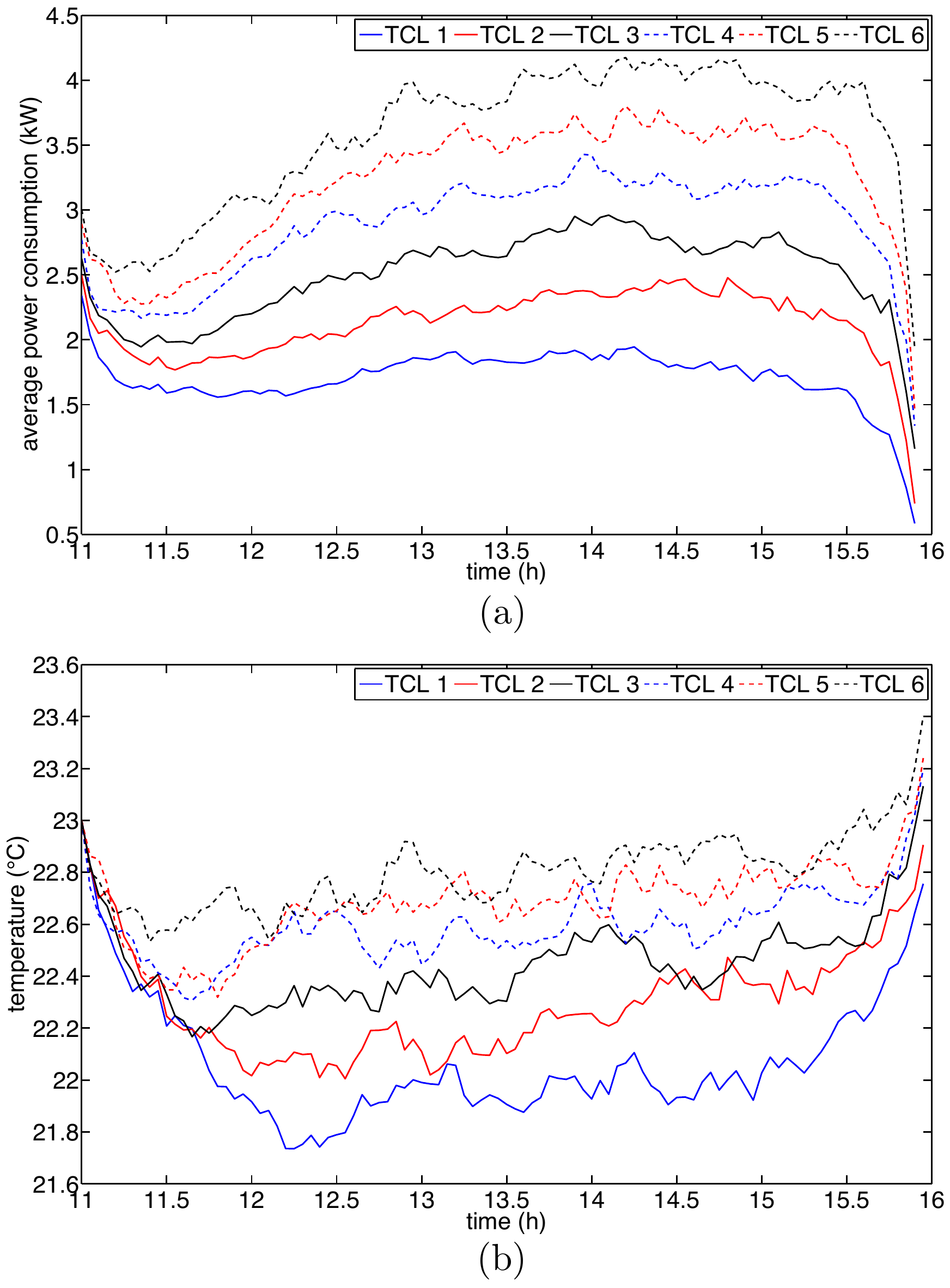}
\caption{(a) Computed optimal control and (b) corresponding state trajectory.
}
 \label{fig:state6}
 \end{center}
\end{figure}

\section{Conclusions and Future Work}

We generalize the path integral-based method for stochastic optimal control to handle a more general class of cost functions, by introducing a new state.
Our generalization makes path integral control more broadly applicable, in particular to applications in engineering, economics and finance.
Furthermore, our memory efficient methodology is applicable to problems where most current techniques fail, because it does not require a global grid for discretizing the domain of the HJB equation.

In the future, we plan to relax the condition on the noise covariance and control input matrices. Another interesting direction would be to develop the path integral method for differential games.

\begin{ack}
The authors would like to thank Professor Lawrence C. Evans for advice on the idea of introducing a new state variable and Professor Duncan S. Callaway for helpful discussions on direct load control.
\end{ack}

\bibliographystyle{alpha}        
\bibliography{general}

\end{document}